\documentclass[12pt,a4paper]{article}

\usepackage[latin1]{inputenc}
\usepackage[T1]{fontenc}
\usepackage[french]{babel} 
\usepackage{amsmath}
\usepackage{amsfonts}
\usepackage{eufrak}
\usepackage{amssymb}
\usepackage{psfig}
\usepackage{epsf}
\usepackage{rotate}
\usepackage[all]{xy}             
\usepackage{xy}
\xyoption{all} 
    
\newcommand{\J}{\mathcal{J}}
\newcommand{\OO}{\mathcal{O}}
\newcommand{\CC}{\mathbb{C}}
 
\newtheorem{defi}{ D\'efinition}[section]

\newtheorem{prop}[defi]{ Proposition}
\newtheorem{coro}[defi]{ Corollaire}

\newtheorem{lemm}[defi]{ Lemme}

\begin{document}
 \setcounter{page}{1}

 \title{Sur la dimension de l'ensemble des points base du fibr\'e d\'eterminant sur $\mathcal{SU}_{C}(r)$}

\author{SCHNEIDER Olivier \thanks{Laboratoire J.-A. Dieudonn\'e
U.M.R. no 6621 du C.N.R.S.
 Universit\'e de Nice - Sophia Antipolis
 Parc Valrose
 06108 Nice Cedex 02
 France.
email : \texttt{oschneid@math.unice.fr}}}
\date{}

 \maketitle
\begin{abstract}

\end{abstract}

 \section{Introduction}
 Soit  $C$  une courbe lisse de genre $g \geq 2$ sur $\CC$. Soit $\J$ la Jacobienne de $C$.
Soit $N$ un fibr\'e ample sur $\J$. Dans ~\cite{R}, M.Raynaud construit \`a partir de $N$, un fibr\'e $F_{N}$ sur $\J$ dont la restriction \`a $C$, not\'ee $E_{N}$, est un fibr\'e semi-stable. On \'etudiera dans un premier temps quelques propri\'et\'es de $F_{N}$ pour en d\'eduire notament la stabilit\'e de $E_{N}$.  \\ 
 Soit $\mathcal{SU}_{C}(r)$ l'espace des modules des fibr\'es vectoriels semi-stables sur $C$, de rang $r$ et de d\'eterminant trivial. On rappelle qu'un fibr\'e $E$ de $\mathcal{SU}_{C}(r)$ est un point base du fibr\'e d\'eterminant si pour tout fibr\'e en droites $L$ de degr\'e $g-1$, $H^{0}(C,E \otimes L) \neq 0$ (voir ~\cite{B2}). Soit $\mathcal{U}_{C}(r,\mu)$ l'espace des modules des fibr\'es vectoriels semi-stables sur $C$, de rang $r$ et de pente $\mu$. Pour tout $\mu \in \mathbb{Z}$, $\mathcal{U}_{C}(r,\mu)$ est, \`a un recouvrement \'etale fini pr\`es, le produit  de $\mathcal{SU}_{C}(r)$ par $\J$. De ce fait, on d\'etermine ici des points bases du fibr\'e d\'eterminant en construisant des fibr\'es $E$ de $\mathcal{U}_{C}(r,g-1)$ tels que pour tout fibr\'e en droites $L$ de degr\'e $0$,  
$$\qquad  \qquad \qquad \qquad \qquad \qquad h^{0}(C,E\otimes L)\geq1. \qquad  \qquad \qquad \qquad \qquad (R)$$ 
Pour cela, on \'etudie le cas $N=\mathcal{O}_{\J}(n\Theta)$, $n\geq 2$, avec $\Theta$ un diviseur de la polarisation principale sur $\J$. On notera dans ce cas $F_{n}$ et $E_{n}$ les fibr\'es associ\'es respectivement sur $\J$ et $C$. Le choix de $N$ n'\'etant pas canonique, les propri\'et\'es d\'emontr\'ees par la suite sur ce fibr\'e seront les v\'erifi\'ees par tout autre fibr\'e $E_{n} \otimes L$, o\`u $L$ est un fibr\'e en droites de degr\'e $0$ sur $C$. On dira que ces fibr\'es sont de "type $E_{n}$". On d\'emontre ici dans la deuxi\`eme partie le r\'esultat suivant : 
\begin{prop} \label{P}
 Soit $n \geq 2$ et $k \geq 1$. Le sous-espace $W$ de $\mathcal{U}_{C}(kn^{g},g-1)$ des fibr\'es stables contenant un fibr\'e de type $E_{n}$ est de codimension $kgn^{2g-1}+1-g$.
\end{prop}
Dans ~\cite{R}, Raynaud d\'emontre que $E_{n}$ est un fibr\'e stable de  rang $n^{g}$, de pente $g/n$, tel que pour tout fibr\'e $L$ de degr\'e $0$, 
$$ h^{0}(C,E_{n}\otimes L)\geq1. $$
De ce fait, pour tout  $r \geq n^{g}$, tous les fibr\'es de $\mathcal{U}_{C}(r,g-1)$ contenant des fibr\'es de type $E_{n}$ ne v\'erifient pas $(R)$ d\`es que $n \geq 2$. Le r\'esultat de la Proposition ~\ref{P} nous donne alors pour $r=kn{g}$, avec $k$ un entier sup\'erieur ou \'egal \`a $1$, une minoration de la codimension du lieu de base du diviseur th\^eta de l'ordre de $\frac{1}{kn} \mathrm{dim}(\mathcal{SU}_{C}(r))$.    

\section {Propri\'et\'es du fibr\'e $F_{N}.$}
  Soit  $C$  une courbe lisse de genre $g \geq 2$ sur $\CC$. Soit $\J$ la Jacobienne de $C$. On identifie canoniquement $\J$ \`a  sa vari\'et\'e ab\'elienne duale gr\^ace \`a la polarisation principale sur $\J$. On note $ \mathcal{P}$ le fibr\'e de Poincar\'e sur $\J \times \J$, trivial sur $\J \times \{0\}$ et sur $\{0\} \times \J$. Soient $p_{1}$ et $p_{2}: \J \times \J \longrightarrow \J$ les deux projections. Soit  $D(\J)$ la cat\'egorie d\'eriv\'ee associ\'ee aux faisceaux coh\'erents sur $\J$ et $\left[ g\right] $  la translation de $g$ degr\'es vers la gauche dans $D(\J)$. Dans ~\cite{Mu}, Mukai d\'efinit la transform\'ee de Fourier-Mukai comme l'application 
$$\begin{array}{rccl} \mathcal{F}: & D(\J) & \longrightarrow &D(\J) \\
                             & M &\longmapsto & Rp_{2,\ast} \left( p_{1}^{\ast}M \otimes \mathcal{P} \right). \end{array} $$ 

Soit $N$ un faisceau ample sur $\J$, on note
$$F_{N}:=\mathcal{F}(N^{-1}) \left[ g\right],$$
ou $F_{D}$ si $N=\mathcal{O}_{\J}(D)$. Soit $\Phi_{N} : \J \longrightarrow \J$ l'isog\'enie qui \`a $x \in \J$ associe $T_{x}^{\ast}N \otimes N^{-1}$, o\`u $T_{x}$ est la translation par $x$ dans $\J$. On note $H(N)$ le noyau de $\Phi_{N}$. Soit $l:=h^{0}(\J,N)$,  $H(N)$ est un sous groupe de $\J$ d'ordre $l^{2}$. On peut d\'efinir sur $H(N)$ une forme bilin\'eaire non d\'eg\'en\'er\'ee \`a valeurs dans $\CC^{\ast}$ (voir ~\cite{M}).
Il existe alors deux sous groupes d'ordre $l$ $K$ et $K^{'}$, totalement isotropes pour cette forme et tels que $H(N)=K \otimes K^{'}$. Soit $\J^{'}:=\J/K^{'}$, c'est une vari\'et\'e ab\'elienne principalement polaris\'ee dont on note $L{'}$ le fibr\'e en droites associ\'e \`a une polarisation principale. On  a une factorisation canonique (voir ~\cite{M})
 $$ \xymatrix{     & \J^{'} \ar[dr]^{\tau} & \\ \J  \ar[rr]^{\phi_{N}} \ar[ur] & & \J} $$
et le r\'esultat suivant : 
 \begin{prop} 
\begin{enumerate}
\item $F_{N}$ est un fibr\'e vectoriel sur $\J$ de rang $l$. \label{m}
\item  $F_{N} \cong \tau_{\ast}L^{'}$ 
\item $\Phi^{\ast}_{N}F_{N} = N \otimes H^{0}(\J,N)^{\ast}$ comme $H(N)$-fibr\'es.
\end{enumerate}
\end{prop}
\begin{description}
\item{\textbf{D\'emonstration}} :  
 Mukai prouve les deux premiers points dans ~\cite{Mu} (p.162). Il \'etablit aussi l'isomorphisme suivant  
$$\qquad  \qquad  \qquad \qquad \qquad \qquad \qquad \Phi^{\ast}_{N}F_{N} \cong N \otimes V    \qquad\qquad \qquad  \qquad \qquad \qquad (\dagger),$$
o\`u $V$ est un espace vectoriel de dimension $l$. On rappelle maintenant quelques r\'esultats dus \`a Mumford (voir ~\cite{M}, p.288--297) : Soit $\mathcal{G}(N)$ l'ensemble des paires $(x,\varphi)$, o\`u $x$ est un \'el\'ement de $H(N)$ et $\varphi$ est un isomorphisme entre $N$ et $T_{x}^{\ast}N$.
$\mathcal{G}(N)$ est un groupe et on a la suite exacte suivante
$$0 \longrightarrow \CC^{\ast} \longrightarrow \mathcal{G}(N) \longrightarrow H(N) \longrightarrow 0.$$

Mumford \'etablit que toute repr\'esentation irreductible de $\mathcal{G}(N)$ sur laquelle $\CC^{\ast}$ agit par homoth\'eties, est isomorphe \`a $\bigoplus_{r} H^{0}(\J,N)$ pour un entier $r$. Du fait de l'isomorphisme $(\dagger)$, $V=H^{0}(C,\Phi^{\ast}F_{N} \otimes N^{-1})$ et comme $H(N)$ agit sur $H^{0}(\J,\Phi^{\ast}_{N}F_{N})$, $V$ est une repr\'esentation irr\'eductible de $\mathcal{G}(N)$ sur laquelle $\CC$ agit par multiplication par l'inverse. $V^{\ast}$ est donc une repr\'esentation de $\mathcal{G}(N)$ sur laquelle $\CC$ agit naturellement ; pour des raisons de dimension, $V^{\ast}$ est isomorphe \`a $H^{0}(\J,N)$, d'o\`u l'\'egalit\'e de $H(N)$-fibr\'es : 
$$\Phi^{\ast}_{N}F_{N} = N \otimes H^{0}(\J,N)^{\ast}.$$ \hfill{$\Box$}
\end{description} 
Soit $E_{N}:=F_{N \mid C}$.  
\begin{coro}
 $E_{N}$ est stable.
\end{coro}

\begin{description}
\item{\textbf{D\'emonstration}} : Soit $C^{''}:= \Phi^{-1}C$. Soit $F$ un sous fibr\'e de $E_{N}$. Alors $$\Phi^{\ast}G=M \otimes W, $$
avec $M$ un sous fibr\'e de $N_{\mid C^{''}}$ et $W$ un sous-espace vectoriel de $H^{0}(\J,N)^{\ast}$. Ceci impose $\mu(F) \leq \mu(E_{N})$. Suppossons maintenant que $F$ est un sous fibr\'e strict de $E_{N}$ tel que $\mu(F)=\mu(E_{N})$, alors
 $$\Phi^{\ast}_{N}F = N_{|C^{''}} \otimes W,$$
avec $W$ un sous-espace vectoriel strict de $H^{0}(\J,N)^{\ast}$ stable sous l'action de  $\mathcal{G}(N)$. Mais ceci est en contradiction avec le fait que $H^{0}(\J,N)^{\ast}$ est une repr\'esentation irr\'eductible de  $\mathcal{G}(N)$. \hfill{$\Box$} \end{description}
Soit maintenant $\mathcal{O}_{\J}(\theta)$, le fibr\'e associ\'e \`a une polarisation principale sur $\J$. Pour tout $x$ dans $\J$ on note $M_{x}:= \Phi_{\mathcal{O}_{\J}(\theta)}(x)$.
\begin{prop} \label{surj}
\begin{enumerate}
\item $\Phi_{N}^{\ast}(\bigoplus_{x\in H(N)} M_{x})=H^{0}(\J,N) \otimes H^{0}(\J,N)^{\ast}\otimes \mathcal{O}_{\J}$ comme $H(N)$-fibr\'e.
\item Si $N$ est engendr\'e par ses sections globales, on a 
$$\bigoplus_{x\in H(N)} M_{x} \longrightarrow F_{N} \longrightarrow 0.$$
\end{enumerate}
\end{prop}
\begin{description}
\item{\textbf{D\'emonstration}} : \begin{enumerate} 
\item Comme $H(N)$ est un groupe commutatif, on a la d\'ecomposition suivante (voir ~\cite{S}): 
$$H^{0}(\J,N) \otimes H^{0}(\J,N)^{\ast}= \bigoplus_{\chi \in Hom(H(N), \CC^{\ast})} V_{\chi},$$
o\`u les $V_{\chi}$ sont des $\CC$-espaces vectoriels de dimension 1 tels que :
$$\forall x \in H(N), \ \forall s \in V_{\chi}, \ x.s=\chi(x)s.$$
De plus, pour tout $x$ dans $H(N)$, $\Phi_{N}^{\ast}(M_{x})$ est un fibr\'e trivial de rang 1 stable sous l'action de $H(N)$, donc
$$\Phi_{N}^{\ast}(M_{x})=V_{\chi_{x}}\otimes \mathcal{O}_{\J},$$
o\`u $\chi_{x}$ est donn\'e par  l'isomorphisme canonique 
$H(N)  \tilde{\longrightarrow}   Hom(H(N), \CC^{\ast}).$ De ce fait, 
$$\Phi_{N}^{\ast}(\bigoplus_{x\in H(N)} L_{x})=H^{0}(\J,N) \otimes H^{0}(\J,N)^{\ast}\otimes \mathcal{O}_{\J}.$$
\item L'application 
$$H^{0}(\J,N) \otimes \mathcal{O}_{\J} \longrightarrow N$$ est surjective. Il en est donc de m\^eme pour $$H^{0}(\J,N) \otimes  H^{0}(\J,N)^{\ast}\otimes \mathcal{O}_{\J}\longrightarrow N \otimes  H^{0}(\J,N)^{\ast}.$$ D'apr\`es ce qui  pr\'ec\`ede, on obtient donc $$\bigoplus_{x\in H(N)} M_{x} \longrightarrow F_{N} \longrightarrow 0.$$ \hfill{$\Box$}
\end{enumerate}
\end{description}

\section{D\'emonstration de le Proposition \ref{P}}
Dans tout ce qui suit, $N:=\mathcal{O}_{\J}(n \Theta)$, o\`u $\Theta$ est un diviseur associ\'e \`a une polarisation principale. On note alors $F_{n}$ et $E_{n}$ les fibr\'es associ\'es. On dira qu'un fibr\'e est de "type $E_{n}$" s'il s'\'ecrit $E_{n} \otimes L$, avec $L$ un fibr\'e en droites de degr\'e $0$ sur $C$.  Le choix de $N$ n'\'etant pas canonique, ces fibr\'es v\'erifient les m\^emes propri\'et\'es que $E_{n}$. Donnons maintenant quelques r\'esultats utiles \`a la d\'emonstration de la Proposition \ref{P}.
\begin{prop}(A.Beauville) \label{a}
Soit  $C$  une courbe lisse de genre $g \geq 2$ sur $\CC$ et $C^{'}$ une autre courbe lisse sur $\CC$. Soit $\pi : C^{'} \longrightarrow C$ un rev\^etement galoisien de groupe de Galois $G$. Pour $M$ un fibr\'e vectoriel stable sur $C^{'}$, on a :
\begin{enumerate}
\item $\pi_{\ast} M$ est semi-stable.
\item Si $\forall g \in G, \ g^{\ast}M \ncong M$, alors $\pi_{\ast}M$ est stable sur $C$. 
\end{enumerate}

\end{prop}

Ce r\'esultat n'\'etant pas pas publi\'e nous en donnons ici une d\'emonstration :
\begin{description}
\item{\textbf{D\'emonstration}} : Remarquons tout d'abord que pour tout faisceau coh\'erent $M$ sur $C^{'}$, on a l'isomorphisme canonique suivant :
$$\pi^{\ast}\pi_{\ast}M \cong \oplus_{g \in G} g^{\ast}M.$$
Soit $M$ un fibr\'e vectoriel stable sur $C^{'}$, comme l'action de $G$ conserve la pente et la stabilit\'e $\pi^{\ast}\pi_{\ast}M$ est somme directe de fibr\'es stables de m\^eme pente. 
\begin{enumerate}
\item Supposons que $\pi_{\ast} M$ contienne un sous fibr\'e $F$ de pente strictement sup\'erieure \`a la pente de $\pi_{\ast}M$. $\pi^{\ast}F$ est alors un sous fibr\'e de $\pi^{\ast}\pi_{\ast}M$ de pente strictement sup\'erieure \`a la pente de $\pi^{\ast}\pi_{\ast}M$, ce qui n'est pas possible.
\item Supposons maintenant que  $\forall g \in G, \ g^{\ast}M \ncong M $ et que $\pi_{\ast}M$ contienne un sous fibr\'e non trivial $F$ tel que $\mu(F)=\mu(\pi_{\ast}M)=\mu$. Alors $\pi^{\ast}F$ est un sous fibr\'e de $\oplus_ {g \in G} g^{\ast}M$ de pente $\mathrm{deg}(\pi)\mu$. Comme les fibr\'es $g^{\ast}M$ sont stables et non isomorphes deux \`a deux, $F=\oplus_{g \in G^{'}} g^{\ast}M$, avec $G^{'}$ un sous ensemble de $G$. Cependant, comme $\pi^{\ast} F$ est $G-$ invariant cela impose $G^{'}=G$ et $F=\pi_{\ast}M$, ce qui prouve la stabilit\'e de $\pi_{\ast}M$.   
\end{enumerate}
\hfill{$\Box$}
\end{description}

\begin{coro} \label{b}
Soit $M$ un fibr\'e vectoriel sur $C^{'}$. Pour $M$ g\'en\'erique, $\pi_{\ast} M $ est stable.
\end{coro}
\begin{description}
\item[D\'emonstration :] remarquons tout d'abord que pour montrer le r\'esultat pour un fibr\'e vectoriel g\'en\'erique de rang $r$, il suffit de le montrer pour un seul fibr\'e vectoriel de rang $r$. Pour cela, on prend comme fibr\'e l'image directe d'un fibr\'e en droites par un morphisme \'etale de degr\'e $r$ sur $C^{'}$. Cela nous am\`ene \`a traiter uniquement le cas des fibr\'es en droites. \\
D'apr\`es la Proposition pr\'ec\'edente, il faut donc  prouver que pour $L$ g\'en\'erique, $$\forall g \in G, \ g^{\ast} L \ncong L. $$
Soit  $g\in G$ et $Z$ la sous vari\'et\'e de $\mathrm{Pic}(C^{'})$ stable par $<g>$. Soit
$$\tau : C^{'}\longrightarrow C^{'}/<g>,$$
alors $Z$ doit donc \^etre contenue dans $\tau^{\ast}\mathrm{Pic}(C^{'}/<g>)$.Mais comme par la formule de Riemann-Hurwitz, le genre de $C^{'}/<g>$ est strictement inf\'erieur au genre de $C^{'}$,
$$\mathrm{dim}(Z) \leq \mathrm{dim}(\tau^{\ast}\mathrm{Pic}(C^{'}/<g>)) < g^{'}=\mathrm{dim}(\mathrm{Pic}(C^{'}).$$
De ce fait pour $L$ suffisament g\'en\'eral, $L$ ne peut \^etre fix\'e par un \'el\'ement de $<g>$.
\hfill{$\Box$}

\end{description}

\begin{lemm}(A.Beauville)  \label{G}
Soit $C$ une courbe lisse de genre $g \geq 2$ sur $\CC$. Soient $L$ un faisceau coh\'erent et $M$ un fibr\'e vectoriel sur $C$ tels que $h^{0}(C,L) \leq h^{1}(C,M)$. Pour $c$ g\'en\'erique dans $\mathrm{Ext}^{1}_{\OO_{C}}(L,M)$, l'application multiplication par $c$
$$H^{0}(C,L) \longrightarrow H^{1}(C,M), $$ est injective.
\end{lemm}

\begin{description}
\item[D\'emonstration :] Comme par dualit\'e $H^{1}(C,M) \cong H^{0}(C,\omega_{C} \otimes M^{\ast})^{\ast}$ et que $\mathrm{Ext}^{1}_{\OO_{C}}(L,M) \cong H^{0}(C, \omega_{C} \otimes L \otimes M^{\ast}) ^{\ast}$, Il revient au m\^eme de d\'emontrer que
$$\xymatrix{ H^{0}(C,L)\otimes H^{0}(C, \omega_{C} \otimes M^{\ast}) \ar[r]^{\ \ \ \ \varphi} &  H^{0}(C, \omega_{C} \otimes L \otimes M^{\ast}) \ar[r]^{\ \ \ \ \  \ \ \ \ \ \ \ c}  & \CC}$$
est s\'eparante \`a gauche pour $c$ g\'en\'erique. \\
Soient $A$, $B$ et $V$ des espaces vectoriels de dimensions respectives $a,b,n$ tels que $a \leq b$. Soit
$$\varphi : A \times B \longrightarrow V,$$ une application bilin\'eaire int\`egre et $c:V \longrightarrow \CC$ g\'en\'erique.
D\'emontrons que $c \circ \varphi$ est s\'eparante \`a gauche. Pour cela, on montre qu'un hyperplan g\'en\'erique de $V$ ne contient aucun des $b-$plans $\varphi(a,b)$. Ces $b-$plans sont param\'etr\'es par une sous-vari\'et\'e $G$ d'une Grassmanienne de dimension inf\'erieure \`a $a-1$. Notons $\breve{P}$ l'espace projectif des hyperplans de $V$ : Soit
$$Z:=\{(L,H) \in G \times \breve{P} \mid L \subset H \}.$$
$Z$ est un fibr\'e sur $G$ de fibre en $L$, l'espace projectif des hyperplans de $V$ contenant $L$, de dimension $n-b-1$ ; de ce fait,
$$\dim Z \leq (a-1)+(n-b-1)=n-2-(b-a) < \dim \breve{P}. $$ \hfill{$\Box$}
\end{description}

\begin{prop} \label{inj}
Soit $n \geq 2$. Soit $\xymatrix{ 0 \ar[r] & E_{n } \ar[r]^{i} &  F \ar[r] & G \ar[r] & 0}$ une extension g\'en\'erique de pente $g-1$  de $E_{n }$ par un faisceau coh\'erent $G$ sur $C$. Alors $i$ est l'unique injection de $E_{n}$ (\`a automorphismes de $E_{n}$ pr\`es) dans $F$ dans les deux cas suivants :
\begin{enumerate}
 \item $G$ fibr\'e vectoriel sur $C$, g\'en\'erique
 \item $G=\oplus_{p\in \mathcal{P}} \CC_{p}$, o\`u $\CC_{p}$ est le faisceau concentr\'e en $p$ et $\mathcal{P}$ est  un ensemble de points de $C$.
 \end{enumerate}
 \end{prop}
\begin{description}
\item{\textbf{D\'emonstration}}: Pour tout $x$ dans $\J_{n}:=\mathrm{Ker}(\Phi_{n \theta})$, on note $L_{x}:=M_{x \mid C}$. On a dans les deux cas \'evoqu\'es dans la Proposition, pour tout $x$ dans $\J_{n}$,
$$h^{1}(C,G \otimes L_{x}^{-1})=0.$$
Mais comme pour tout fibr\'e en droites $L$ de degr\'e $0$,
$$\chi(G \otimes L)=- \chi(E_{n} \otimes L) >0,$$
ceci implique
$$h^{0}(C,G \otimes L_{x}^{-1}) \leq h^{1}(C,E_{n} \otimes L_{x}^{-1}).$$
Par g\'en\'ericit\'e de l'extension $F$ dans $\mathrm{Ext}^{1}_{\OO_{C}}(G,E_{n }) \cong H^{0}(C,\omega_{C} \otimes (E_{n} \otimes L_{x}^{-1})^{\ast} \otimes (G \otimes L_{x}^{-1}))^{\ast}$, on peut donc appliquer le Lemme \ref{G} pour obtenir que pour tout $x$ dans $\J_{n}$, l'application $$H^{0}(C,G\otimes L_{x}^{-1}) \longrightarrow H^{1}(C,E_{n} \otimes L_{x}^{-1})$$ est injective. De ce fait, en \'ecrivant la suite exacte longue d'homologie, on obtient que   $$H^{0}(C,E_{n}\otimes  L_{x}^{-1}) \cong  H^{0}(C,F\otimes  L_{x}^{-1}).$$
En cons\'equence, le diagramme suivant est commutatif :
$$ \xymatrix{ \bigoplus_{x \in \J_{n}} \mathrm{Hom}(L_{x},E_{n}) \otimes L_{x} \ar[r]^{ \quad \quad \quad \quad \quad \psi}  \ar[d]^{\wr}&  E_{n} \ar@{^{(}->}[d]^{i}  \\ \bigoplus_{x \in \J_{n}} \mathrm{Hom}(L_{x},F) \otimes L_{x}  \ar[r]^{ \quad \quad \quad \quad \quad \varphi}  &  F. }$$
Mais comme pour $n \geq 2$, $\mathcal{O}_{\J}(n \theta)$ est engendr\'e par ses sections globales, on a par la Proposition \ref{surj},
$$\bigoplus_{x\in \J_{n}} M_{x} \longrightarrow F_{n} \longrightarrow 0.$$
 De ce fait, $\psi$ est surjective et l'image de $i$ est \'egale \`a $\mathrm{Im}(\varphi)$.    \hfill{$\Box$}
\end{description}

\begin{description}
\item{\textbf{D\'emonstration de la Proposition \ref{P}} } :
Consid\'erons une suite exacte
$$\xymatrix{ 0 \ar[r] & E_{n } \ar[r] &  F \ar[r] & G \ar[r] & 0},$$ pour les deux cas suivants :

\begin{enumerate}
\item $G=\oplus_{p\in \mathcal{P}} \CC_{p}$, $\CC_{p}$ \'etant le faisceau concentr\'e en $p$ et $\mathcal{P}$  un ensemble de points de $C$ de cardinal $(g-1)n^{g}-gn^{g-1}$.
\item $G$ fibr\'e stable de degr\'e $k(g-1)n^{g}-gn^{g-1}$ et de rang $(k-1)n^{g}$ avec $k \geq 2$.
\end{enumerate}
Cela impose $\mu(F)=g-1$. D\'emontrons que les fibr\'es ainsi construits sont g\'en\'eriquement stables. Pour cela, on consid\`ere  $\pi:C^{'} \longrightarrow C$ la restriction de $\tau:\J^{'} \longrightarrow \J$ \`a $C^{'}:=\tau^{-1}C$ et $K:=\mathrm{Ker}(\pi)$. D'apr\`es la Proposition \ref{m}, $E_{n }=\pi_{\ast}L^{'}$, o\`u $L^{'}$ est la restriction  \`a $C^{'}$ d'un fibr\'e associ\'e \`a une polarisation principale de $\J^{'}$. On va maintenant d\'emontrer dans les deux cas, l'existence d'un fibr\'e stable $F$, image directe par $\pi$ d'une extension de $L^{'}$ :
\begin{enumerate}
\item Si $k=1$, consid\'erons le fibr\'e en droites $F^{'}:=L^{'}(\sum_{p\in \mathcal{P}}p)$, o\`u $\mathcal{P}$  est un ensemble de points de $C^{'}$ de cardinal $(g-1)n^{g}-gn^{g-1}$. Alors o\`u bien pour tout $\sigma \in K$, $\sigma^{\ast}F^{'} \ncong F^{'}$ et alors $F^{'}$ v\'erifie les conditions de  la Proposition \ref{a} et de ce fait $F:=\pi_{\ast}F^{'}$ est stable. O\`u bien il existe $\sigma \in K$ tel que $\sigma^{\ast}F^{'} \cong F^{'}$. Dans ce cas, soit $q_{1} \in \mathcal{P}$, alors il existe $q_{2} \in  C^{'}- \mathcal{P}$ tel que pour tout $\sigma \in K$,
$$\sigma^{\ast} L^{'}( \sum_{p\in \mathcal{P}} p) \otimes   L^{'}( \sum_{p\in \mathcal{P}}p )^{-1} \ncong \OO_{C^{'}}(q_{1}-\sigma(q_{1})+\sigma(q_{2})-q_{2}).$$
De ce fait, pour tout $\sigma \in K$, $\sigma^{\ast}L^{'}(q_{2}-q_{1}+\sum_{p\in \mathcal{P}}p) \ncong L^{'}(q_{2}-q_{1}+\sum_{p\in \mathcal{P}}p) $ et donc par la Proposition \ref{a} son image directe par $\pi$ sera stable.
\item  Consid\'erons maintenant un fibr\'e $F^{'}$ extension de $L^{'}$ par un fibr\'e vectoriel $G^{'}$. $\det(F^{'})=L^{'} \otimes \det(G^{'})$, donc pour $G^{'}$ g\'en\'erique, pour tout $\sigma \in K$, $\sigma^{\ast}\det(F^{'}) \ncong \det(F^{'})$. Cette propri\'et\'e est aussi v\'erifi\'ee par $F^{'}$ et par la Proposition \ref{a} ceci implique que dans ce cas $\pi_{\ast} F^{'}$ est stable.
\end{enumerate}

Effectuons maintenant le d\'ecompte des dimensions : par la Proposition \ref{inj}, il suffit d'ajouter la dimension de l'espace des fibr\'es de type $E_{n}$, c'est \`a dire $g$, \`a la dimension de l'espace dans lesquel on choisit $G$ plus $\mathrm{dim}(\mathrm{Ext}^{1}_{\OO_{C}}(G,E_{n}))-\mathrm{dim}(\mathrm{Aut}(G))$ :
\begin{enumerate}
\item Si $G=\oplus_{p\in \mathcal{P}} \CC_{p}$, $\CC_{p}$, alors
$$\mathrm{dim}(\mathrm{Ext}^{1}_{\OO_{C}}(\oplus_{p\in \mathcal{P}} \CC_{p}, E_{n}))=\mathrm{rg}(E_{n }) \left((g-1)n^{g}-gn^{g-1}\right)=n^{g}\left((g-1)n^{g}-gn^{g-1}\right).$$
Comme $\mathrm{dim}(\mathrm{Aut}(G))=\mathrm{Card}(\mathcal{P})$ et que la dimension de l'espace dans lequel on choisit le fibr\'e $G$ est aussi $\mathrm{Card}(\mathcal{P})$,  on a :
$$ \mathrm{dim}(W) = g+n^{g}\left((g-1)n^{g}-gn^{g-1}\right). $$

\item Si $k \geq 2$,  $\mathrm{dim}(\mathrm{Ext}^{1}_{\OO_{C}}(G,E_{n}))=h^{1}(C,G^{\ast} \otimes E_{n})$. Or pour $G$ g\'en\'erique et stable de pente sup\'erieure \`a la pente de $E_{n}$, on a $h^{0}(C,G^{\ast} \otimes E_{n\theta})=0$. De ce fait, on obtient par Riemann-Roch,
$$\begin{array}{rcl} h^{1}(C,G^{\ast} \otimes E_{n}) & = & \mathrm{rg}(E_{n}) \mathrm{deg}(G)-\mathrm{deg}(E_{n}) \mathrm{rg}(G)+\mathrm{rg}(E_{n}) \mathrm{rg}(G)(g-1)\\
& = & n^{g}(k(g-1)n^{g}-gn^{g-1})-(k-1)n^{g}gn^{g-1}+n^{g}(k-1)n^{g}(g-1)\\
& = & \left((kn^{g})^{2}(g-1)+1\right)-\left(((k-1)n^{g})^{2}(g-1)+1  \right)-gkn^{2g-1}. \end{array}$$
et comme $G$ est g\'en\'erique dans $U_{s}((k-1)n^{g},\frac{1}{k-1}(g-1-\frac{g}{n}))$, $\mathrm{dim}(\mathrm{Aut}(G))=1$ et donc
$$\mathrm{dim}(W)=g+(kn^{g})^{2}(g-1)-\left(gkn^{2g-1}\right).$$

\end{enumerate}
\hfill{$\Box$}
\end{description}


\end{document}